\title{\bf{New examples of canonical covers of degree 3}}
\author{NGUYEN BIN}
\date{}

\newcommand{\Addresses}{{
		\bigskip
		\footnotesize
		\text{Department of Mathematics and Statistics,}\par\nopagebreak
		\text{Quy Nhon University,}\par\nopagebreak
		\text{170 An Duong Vuong, Quy Nhon, Vietnam.}\par\nopagebreak
		\textit{E-mail address}: \texttt{nguyenbin@qnu.edu.vn}		
	}}

\newcommand{\Addressesb}{{
		\bigskip
		\footnotesize
		\text{Center for Mathematical Analysis, Geometry and Dynamical Systems,}\par\nopagebreak
		\text{Departamento de Matem\'{a}tica,}\par\nopagebreak
		\text{Instituto Superior T\'{e}cnico, Universidade de Lisboa}\par\nopagebreak
		\text{Av. Rovisco Pais, 1049-001 Lisboa, Portugal.}	
}}

\newcommand{\Addressesc}{{
		\bigskip
		\footnotesize
		\text{Mathematics Division,}\par\nopagebreak
		\text{National Center for Theoretical Sciences,}\par\nopagebreak
		\text{Taipei 10617, Taiwan.}	
}}

\documentclass[10pt]{article}
\usepackage{amsmath, amsthm, amssymb, mathrsfs, amscd, amsthm, amscd,amsfonts,calligra,mathrsfs,nicefrac}
\usepackage{adjustbox}
\usepackage{graphicx}

\usepackage{indentfirst,url,xypic}
\usepackage[all]{xy}
\usepackage{url}
\usepackage{tikz}
\usepackage[colorlinks,plainpages]{hyperref}
\hypersetup{
	colorlinks=true,
	linkcolor=blue,
	filecolor=magenta,      
	urlcolor=cyan,
}

\DeclareMathOperator{\degree}{deg}

\newtheorem{Definition}{Definition}[section]
\newtheorem{Theorem}[Definition]{Theorem}
\newtheorem{Remark}[Definition]{Remark}
\newtheorem{Proposition}[Definition]{Proposition }

\newcommand{\MSC}{\textbf{Mathematics Subject Classification (2010):}}

\begin{document}
\maketitle
\begin{abstract} This paper presents new examples of projective surfaces of general type over $\mathbb{C}$ with canonical map of degree $ 3 $ onto a surface of general type. Very few examples are known of such surfaces and some of the examples in this paper present the new feature of having the canonical map not a morphism (i.e. the canonical linear system with base points).\\
	
	\noindent
	\MSC{14J29}.
	
\end{abstract}

\section{Introduction}
Let $ X $ be a minimal smooth complex surface of general type (see \cite{MR1406314} or \cite{MR2030225}) and denote by $ \varphi_{\left| K_X\right| }:X \dashrightarrow \mathbb{P}^{p_g\left( X\right)-1}$ the canonical map of $ X $, where $ p_g\left( X\right) = \dim H^0\left( X, K_{X}\right) $ is the geometric genus and $ K_{X} $ is the canonical divisor of $ X $. A classical result of A. Beauville \cite[\rm Proposition 4.1 and its proof]{MR553705} says that if the image of $ \varphi_{\left | K_X\right| } $ is a surface of general type, the degree $ d $ of the canonical map of $ X $ is less than or equal to $ 9 $. These surfaces are called canonical covers. The problem of constructing canonical covers of degree $ d >1 $ has a long history. At first it was erroneously believed that canonical covers of degree $ d > 1$ did not exist, but around 1979 a first example of a canonical cover of degree $ d = 2 $ was discovered independently by F. Catanese, G. Van der Geer, D. Zagier and  A. Beauville \cite{MR0485704}, \cite{MR620679}. Nowadays, infinitely many canonical covers of degree $ d=2 $ have been constructed by A. Beauville \cite{MR927956}, by C. Ciliberto, R. Pardini and F. Tovena \cite{MR1822412}, \cite{MR1960801} and the author \cite{2019arXiv190713584B}, etc.\\

Only a few examples of canonical covers of degree $ d $ greater than or equal to $ 3 $ have been found so far. One example with $ d = 5 $ has been constructed by S.L. Tan \cite{MR1141782} and R. Pardini \cite{MR1103913} independently. One example with $ d = 3 $, $ p_g = 5 $ and $ K^2 = 27 $ has been found by R. Pardini \cite{MR1103913} in 1991 and a series of examples with $ d = 3 $, $p_g = 8-n $ and $ K^2 = 54- 9n$ for some $ 0 \le n \le 4 $ have been found by S.L. Tan \cite{MR1141782} in 1992. Later, W. Barth \cite[unpublished]{Barth2000} showed that the quintic surface constructed by G. Van der Geer and D. Zagier \cite{MR0485704} in 1977 is also a canonical cover of degree $ 3 $. This quintic surface has $ p_g = 4 $, $ K^2 = 15 $. In 2016, C. Rito constructed another quintic surface with $ p_g = 4 $, $ K^2 = 15 $ \cite{MR3464003}. All the these known canonical covers of degree $ 3 $ possess base point free canonical system.\\

In this paper, we construct several new canonical covers of degree $ 3 $. Some of these new surfaces have base points for the canonical system. The results of this paper are the following theorems:

\begin{Theorem}\label{THE MAIN THEOREM WITH D=3 PART 1}
	There exist minimal surfaces of general type $ X $ whose canonical map is of degree 3 onto a surface of general type $ \Sigma $, satisfying the following
	$$
	\begin{tabular}{|c| c| c| c| c|}
	\hline
	$ K_X^2 $ &$ p_g\left( X\right) $ &$ q\left( X\right) $ & $ \degree\left( \Sigma\right)  $& Base points of $\left| K_X \right| $\\
	\hline
	$ 30 $&$ 5 $&$ 0 $&$ 8 $& $ 6 $ \\
	\hline
	$ 28 $&$ 5 $&$ 0 $&$ 8 $& $ 4 $ \\
	\hline
	$ 26 $&$ 5 $&$ 0 $&$ 8 $& $ 2 $ \\
	\hline
	$ 24 $&$ 5 $&$ 0 $&$ 8 $& $ 0 $ \\
	\hline
	$ 21 $&$ 4 $&$ 0 $&$ 5 $& $ 6 $ \\
	\hline
	$ 19 $&$ 4 $&$ 0 $&$ 5 $& $ 4 $ \\
	\hline
	$ 17 $&$ 4 $&$ 0 $&$ 5 $& $ 2 $ \\
	\hline
	\end{tabular} 
	$$
	Furthermore, the image $ \Sigma $ of the canonical map is either a quintic surface in $ \mathbb{P}^3 $ ($ p_g\left( X\right) =4  $) or a complete intersection of a quadric and a quartic in $ \mathbb{P}^4 $ ($ p_g\left( X\right) =5  $).
\end{Theorem}

\begin{Theorem}\label{THE MAIN THEOREM WITH D=3 PART 2}
	There exists a minimal surface of general type $ X $ whose canonical map is of degree 3 onto a surface of general type $ \Sigma $, satisfying the following
	$$
	\begin{tabular}{|c| c| c| c| c|}
	\hline
	$ K_X^2 $ &$ p_g\left( X\right) $ &$ q\left( X\right) $ & $ \degree\left( \Sigma\right)  $& Base points of $\left| K_X \right| $\\
	\hline
	$ 35 $&$ 6 $&$ 0 $&$ 11 $& $ 2 $ \\
	\hline
	\end{tabular} 
	$$		   
\end{Theorem}

\vskip 0.5cm
The idea of our constructions is the following: we construct the surfaces in the first rows of Theorem \ref{THE MAIN THEOREM WITH D=3 PART 1} and \ref{THE MAIN THEOREM WITH D=3 PART 2} by taking specific $ \mathbb{Z}_{3}^2 $-covers $ X $ of a blow-up of $ \mathbb{P}^2$. In particular, we take the covers to be in such a way that one of the $ \mathbb{Z}_{3} $-quotients $ X_1 $ of $X$ is a surface of general type whose only singularities are of type $ A_2 $ or of type $ \frac{1}{3}\left( 1,1\right)  $. Moreover, we require that: 
\begin{enumerate}
	\item $ p_g\left( \overline{X}_1\right) = p_g\left( X\right)  $,
	\item the canonical map of $ \overline{X}_1 $ is birational,
	\item the canonical map of $ X $ essentially factors through the quotient map $ \xymatrix{X \ar[r] & X_1} $,	
\end{enumerate}
\noindent
where $ \xymatrix{\overline{r}: \overline{X}_1 \ar[r] & X_1} $ is a partial resolution of $ X_1 $ resolving only the singularities of type $ \frac{1}{3}\left(1,1 \right)  $. More precisely, we construct $ \mathbb{Z}_{3}^2 $-covers such that the following diagram commutes:
$$
\xymatrix{X \ar[0,3]^{\mathbb{Z}_3^2} \ar[1,1]^{\mathbb{Z}_3} \ar@{.>}[3,1]_{\varphi_{\left| K_X \right| }}&&& Bl\left( \mathbb{P}^2\right) \\
	&X_1 \ar[-1,2]^{\mathbb{Z}_3} &&\\
	&\overline{X}_1 \ar@{.>}[1,0]^{\varphi_{\left| K_{\overline{X}_1} \right| }} \ar[-1,0]^{\overline{r}}&&\\
	&\hskip 0.5cm \Sigma \subseteq \mathbb{P}^{p_g -1}&&} 	
$$
and such that the canonical map $ \varphi_{\left| K_{\overline{X}_1} \right| } $ is birational. In our constructions, the building data of $ \mathbb{Z}_3^2 $-covers is taken so that $ q\left( X\right)  = 0$ and that the quotient surface $ X_1 $ contains a $ \mathbb{Z}_3 $-set (see Definition \ref{Z_3-set}) of $ n $ singular points of type $ A_2 $ and $ m $ singular points of type $ \frac{1}{3}\left( 1,1\right)  $ with $ 2n + m =6\chi\left( \mathcal{O}_{\overline{X}_1} \right)  $. The $ \mathbb{Z}_3 $-cover $ \xymatrix{X \ar[r] & X_1} $ is only ramified on this $ \mathbb{Z}_3 $-set of singularities. It follows from Proposition \ref{singularities on Z_3 quotient} that $ p_g\left( \overline{X}_1\right) = p_g\left( X\right)  $. The main difference between our constructions and the constructions of Tan and Pardini is the singularities of type $ \frac{1}{3}\left( 1,1\right)  $ which are used in the branch locus of $ \mathbb{Z}_3 $-covers.\\

\noindent
To obtain the other surfaces described in Theorem \ref{THE MAIN THEOREM WITH D=3 PART 1}, we modify the $ \mathbb{Z}_{3}^2 $-cover which we obtain the surface in the first row in Theorem \ref{THE MAIN THEOREM WITH D=3 PART 1}. We impose ordinary triple points on the branch locus of this $ \mathbb{Z}_{3}^2 $-cover and resolve the singularities in such a way that the three conditions above hold. The base points come from the singularities of type $ \frac{1}{3}\left( 1,1\right)  $ of the quotient surfaces.\\

{\bf{Notation and conventions}:} Throughout this paper all surfaces are projective algebraic over the complex numbers. Linear equivalence of divisors is denoted by $ \equiv $. For a finite abelian group $ G $, a character $ \chi $ of $ G $ is a homomorphism from $ G $ to $ \mathbb{C}^{*} $, the multiplicative group of the nonzero complex numbers. The additive group of integers modulo $ 3 $ is denoted by $ \mathbb{Z}_3 $. The rest of the notation is standard in algebraic geometry. 

\section*{Acknowledgments}
The author is deeply indebted to Margarida Mendes Lopes for all her help. Many thanks are due to the anonymous referees for their careful reading and detailed comments with many helpful suggestions. The author was partially supported by Funda\c{c}\~{a}o para a Ci\^{e}ncia e Tecnologia (FCT), Portugal through the program Lisbon Mathematics PhD (LisMath) of the University of Lisbon, scholarship FCT - PD/BD/113632/2015 and project UID/MAT/04459/2019 of CAMGSD. This paper was finished during the author's postdoctoral fellowship at the National Center for Theoretical Sciences (NCTS), Taiwan, under the grant number MOST 109-2119-M-002-014. The author would like to thank NCTS for the financial support and kind hospitality.

\section{$ \mathbb{Z}_{3}$ and $ \mathbb{Z}_{3}^2 $-coverings}
In this section, we recall some facts on abelian covers, in a form which is convenient for our later constructions, in particularly on $ \mathbb{Z}_{3}$-covers and $ \mathbb{Z}_{3}^2$-covers.\\

An abelian cover of a smooth surface $ Y $ with a finite abelian group $ G $ is a finite map $ \xymatrix{f: X \ar[r] & Y} $ together with a faithful action of $ G $ on $ X $ such that $ f $ exhibits $ Y $ as the quotient of $ X $ via $ G $.\\

The construction of abelian covers was studied by R. Pardini in \cite{MR1103912}. An abelian cover $ \xymatrix{f: X \ar[r] & Y} $ with $ X $ normal and $ Y $ smooth can be constructed from a collection of line bundles $ L_{\chi} $ and effective divisors $ D_{\left( H, \psi\right) } $ on the target surface $ Y $. The line bundles $ L_{\chi} $ are labelled by non-trivial characters $ \chi $ of the group $ G $ and the effective divisors $ D_{\left( H, \psi\right) } $ are labelled by pairs $ \left( H, \psi\right) $, where the first index $ H $ is a cyclic subgroup of $ G $ and the second index $ \psi $ is a generator of the group of characters $ H^{*} $.\\

Let $ \xymatrix{f: X \ar[r] & Y} $ be an abelian cover with $ X $ normal and $ Y $ smooth. The set of data $ \left\lbrace L_{\chi}, D_{\left( H, \psi\right) } \right\rbrace_{\chi,\left( H, \psi\right)}$ satisfies 
\begin{align}\label{the relation of an abelian cover}
L_{\chi} + L_{\chi^{'}} \equiv L_{\chi\chi^{'}} + \sum_{H}{\sum_{\psi \in S_H}{\varepsilon}_{\chi, \chi^{'}}^{H,\psi}D_{\left( H, \psi\right) }}, 
\end{align}
for all characters $ \chi $ and $ \chi^{'} $ of $ G $. In the above formula, the first sum runs over all non-trivial cyclic subgroups $ H $ of $ G $, $ S_H $ is the set of generating elements of $ H^{*} $, and the integer number $ \varepsilon_{\chi, \chi^{'}}^{H,\psi} $ is defined as follows:
\begin{align*}
\varepsilon_{\chi, \chi^{'}}^{H,\psi}:=
\left\{ \begin{array}{l}
0 \hskip 0.5cm \text{if } \iota_{\chi} + \iota_{\chi^{'}} < ord\left( H\right) \\
1 \hskip 0.5cm \text{else}
\end{array}  \right.  
\end{align*}
where $ \iota_{\chi}, \iota_{\chi^{'}} \in \left\lbrace 1,2,\ldots, ord\left( H\right)-1\right\rbrace $ are integer numbers such that
\begin{align*}
\chi|_{H} = \psi^{\iota_{\chi}}, \hskip 0.5cm \chi^{'}|_{H} = \psi^{\iota_{\chi^{'}}}.
\end{align*}
\noindent
The set of data $ \left\lbrace L_{\chi}, D_{\left( H, \psi\right) } \right\rbrace_{\chi,\left( H, \psi\right)}$ is called the building data of the abelian cover $ \xymatrix{f: X \ar[r] & Y} $ if $ X $ is normal. By \cite[Theorem 2.1]{MR1103912}, to any set of data $ \left\lbrace L_{\chi}, D_{\left( H, \psi\right) } \right\rbrace_{\chi,\left( H, \psi\right)}$ satisfying $ (\ref{the relation of an abelian cover}) $ we can associate an abelian cover $ \xymatrix{f: X \ar[r] & Y} $ in a natural way.\\

If $ X $ is normal and $ Y $ is smooth, the action of $ G $ induces the following $ G $-equivariant decomposition
\begin{align*}
H^{0}\left( X, K_X\right) = \bigoplus_{\chi \in G^{*}}{H^{0}\left( Y, K_{Y} +L_{\chi^{-1}}\right)} 
\end{align*}
\noindent
where $ G^{*} $ is the character group of $ G $ and the group $ G $ acts on $ H^{0}\left( Y, K_{Y} +L_{\chi^{-1}}\right) $ via the character $ \chi $ (see \cite[ \rm Proposition 4.1c]{MR1103912}). Let $ H $ be a subgroup of $ G $ and let $ H^\perp $ denote the kernel of the restriction map $ \xymatrix{ G^{*} \ar[r]&H^{*}} $. If $ h^0\left( K_{Y} +L_{ \chi} \right) = 0 $ for all $ \chi \notin H^\perp $, the subgroup $ H $ acts trivially on $ H^{0}\left( X, K_X\right) $. So the canonical map of $ X $ essentially factorizes via the quotient map $\xymatrix{X \ar[r]&  X/H} $ (see e.g. \cite[\rm Example 2.1]{MR1103913}). Summarizing, we can now state the following proposition:
\begin{Proposition}\label{Factorization of the canonical map through a quotient surface}
	Let $ \xymatrix{f: X \ar[r] & Y} $ be an abelian $ G $-cover of a smooth surface $ Y $ with $ X $ normal. If there exists a subgroup $ H $ of $ G $ such that $ h^0\left( K_{Y} +L_{ \chi} \right) = 0 $ for all $ \chi \notin H^\perp $, birationally the canonical map of $ X $ is the composition of the quotient map $\xymatrix{X \ar[r]&  X/H} $ with the canonical map of $\overline{X/H}$, where the surface $\overline{X/H}$ is a partial resolution with only canonical singularities of $ X/H $.
\end{Proposition} 
\subsection{$ \mathbb{Z}_{3}^2 $-coverings}
Let $ H_{i_1i_2} $ denote the non-trivial cyclic subgroup generated by $ \left( i_1,i_2\right) $ of $ \mathbb{Z}_{3}^2$ for all $ \left( i_1,i_2\right) \in \mathbb{Z}_{3}^2 \setminus \left( 0,0\right) $ and denote by $ \chi_{j_1j_2} $ the character of $ \mathbb{Z}_{3}^2 $ defined by
\begin{align*}
\chi_{j_1j_2}\left( a_1,a_2\right): =  e^{\left( \frac{2\pi}{3} a_1j_1\right) \sqrt{-1}}e^{\left( \frac{2\pi}{3} a_2j_2\right) \sqrt{-1}}
\end{align*}
for all $ j_1,j_2,a_1,a_2 \in \mathbb{Z}_3 $. The group $ \mathbb{Z}_{3}^2 $ consists of four non-trivial cyclic subgroups, namely: 
\begin{align*}
H_{01} = H_{02}, H_{10} = H_{20}, H_{11} = H_{22}, H_{12} = H_{21}. 
\end{align*}
\noindent
Each of these subgroup has two non-trivial characters, namely:
\begin{align*}
H_{01}^{*}&=\left\lbrace \chi_{01}\mid_{H_{01}}, \chi_{02}\mid_{H_{02}}, \textbf{1} \right\rbrace,\\
H_{10}^{*}&=\left\lbrace \chi_{10}\mid_{H_{10}}, \chi_{20}\mid_{H_{20}}, \textbf{1} \right\rbrace,\\
H_{11}^{*}&=\left\lbrace \chi_{11}\mid_{H_{11}}, \chi_{22}\mid_{H_{22}}, \textbf{1} \right\rbrace,\\
H_{12}^{*}&=\left\lbrace \chi_{12}\mid_{H_{12}}, \chi_{21}\mid_{H_{21}}, \textbf{1} \right\rbrace, 
\end{align*}
\noindent
where $ \chi_{j_1j_2}\mid_{H_{j_1j_2}} $ is the restriction of the character $ \chi_{j_1j_2} $ on the subgroup $ H_{i_1i_2} $ and $ \textbf{1} $ is the trivial character.\\

\noindent
For the sake of simplicity, from now on we use the notations $ L_{j_1j_2} $ instead of $ L_{\chi_{j_1j_2}}$ and $ D_{j_1j_2} $ instead of $ D_{\left( H_{j_1j_2},\chi_{j_1j_2}\mid_{H_{j_1j_2}}\right) }$ for all $ j_1,j_2,a_1,a_2 \in \mathbb{Z}_3 $.\\     

From \cite[\rm Theorem 2.1 and Proposition 2.1]{MR1103912} we have
\begin{Proposition}\label{Construction of bitriple cover}
	Given $ Y $ a smooth projective surface, let $ L_{\chi} $ be divisors of $ Y $ such that $ L_{\chi} \not\equiv \mathcal{O}_Y $ for all non-trivial characters $ \chi $ of $ \mathbb{Z}_{3}^2 $ and let $ D_{\sigma} $ be effective divisors of $ Y $  for all $ \sigma \in \mathbb{Z}_3^2\setminus\left\lbrace \left( 0,0\right)  \right\rbrace  $ such that the total branch divisor $ B:=\sum\limits_{\sigma \ne \left( 0,0\right) }{D_{\sigma}} $ is reduced. If $ \left\lbrace L_{\chi}, D_{\sigma} \right\rbrace_{\chi,\sigma}$ is the building data of a $ \mathbb{Z}_{3}^2$-cover $ \xymatrix{f:X \ar[r]& Y} $, then
	$$
	\begin{adjustbox}{max width=\textwidth}
	\begin{tabular}{l l r l r l r l r l r l r l r l r}
	\color{blue}$ 3L_{10} $&\color{blue}$ \equiv $&&$   $&&\color{blue}$ $ &\color{blue} $D_{10} $&\color{blue}$ +$ &\color{blue} $2D_{20} $&\color{blue}$ +$&\color{blue}$2D_{11} $&\color{blue}$ +$&\color{blue}$D_{22} $&\color{blue}$ +$&\color{blue}$2D_{12} $&\color{blue}$ +$&\color{blue}$D_{21} $\\
	\color{blue}$ 3L_{01} $&\color{blue}$ \equiv $ &\color{blue} $D_{01}$&\color{blue}$ +$ &\color{blue} $2D_{02}  $&$  $&&$  $&&\color{blue}$ +$&\color{blue}$2D_{11} $&\color{blue}$ +$&\color{blue}$D_{22} $&\color{blue}$ +$&\color{blue}$D_{12} $&\color{blue}$ +$&\color{blue}$2D_{21} $ \\
	$ 3L_{20} $&$ \equiv $&&$ $&&$$ & $2D_{10} $&$ +$ & $D_{20} $&$ +$&$D_{11} $&$ +$&$2D_{22} $&$ +$&$D_{12} $&$ +$&$2D_{21} $ \\
	$ 3L_{02} $&$ \equiv$ & $2D_{01} $&$ +$ & $D_{02} $&$  $&&$  $&&$ +$&$D_{11} $&$ +$&$2D_{22} $&$ +$&$2D_{12} $&$ +$&$D_{21} $ \\
	$ 3L_{11} $&$ \equiv $ & $D_{01}$&$ +$ & $2D_{02}  $&$ +$ & $D_{10} $&$ +$ & $2D_{20} $&$ +$ & $D_{11} $&$ +$&$2D_{22} $&&$  $&&$  $ \\
	$ 3L_{22} $&$ \equiv $ & $2D_{01}$&$ +$ & $D_{02}  $&$ +$ & $2D_{10} $&$ +$ & $D_{20} $&$ +$ & $2D_{11} $&$ +$&$D_{22} $&&$  $&&$  $ \\
	$ 3L_{12} $&$ \equiv $ & $2D_{01}$&$ +$ & $D_{02}  $&$ +$ & $D_{10} $&$ +$ & $2D_{20} $&$  $&&$  $&&$ +$&$D_{12} $&$ +$&$2D_{21} $ \\
	$ 3L_{21} $&$ \equiv $ & $D_{01}$&$ +$ & $2D_{02}  $&$ +$ & $2D_{10} $&$ +$ & $D_{20} $&$  $&&$  $&&$ +$&$2D_{12} $&$ +$&$D_{21} $.
\end{tabular}
\end{adjustbox}
$$
Conversely, if $ \left\lbrace L_{\chi}, D_{\sigma} \right\rbrace_{\chi,\sigma}$ satisfies the above conditions and the surface $ Y $ has no non-trivial $ 3 $-torsion, then we can associate a $ \mathbb{Z}_{3}^2$-cover $ \xymatrix{f:X \ar[r]& Y} $ in a natural way.
\end{Proposition}
\noindent
In the above proposition, the last statement follows from the fact that if the surface $Y$ has no non-trivial $ 3 $-torsion, then $3L\equiv D$ uniquely defines $L$.\\

\noindent
To obtain the examples described in Theorems \ref{THE MAIN THEOREM WITH D=3 PART 1} and \ref{THE MAIN THEOREM WITH D=3 PART 2}, we do not consider the most general $  \mathbb{Z}_{3}^2 $-covers. We take $ D_{\sigma} $ to be empty for some $ \sigma \in \mathbb{Z}_3^2 $.\\

From \cite[\rm Proposition 3.1]{MR1103912} we have:
\begin{Proposition} \label{Smoothness condition of bitriple cover}
Let $ Y $ be a smooth surface and let $ \xymatrix{f: X \ar[r]& Y} $ be a $  \mathbb{Z}_{3}^2 $-cover with the building data $ \left\lbrace L_{\chi}, D_{\sigma} \right\rbrace_{\chi,\sigma}$ such that $ X $ is normal. Then the surface $ X $ is smooth above a point $ y \in Y$ if and only if one of the following conditions holds:
\begin{enumerate}
	\item the point $ y $ is not a branch point of $ f $.
	\item the point $ y $ belongs only to one component $ \Delta $ of $ \cup_{\sigma}D_{\sigma} $ and $ y $ is a smooth point of $ \cup_{\sigma}D_{\sigma} $.
	\item the point $ y $ belongs only to $ D_{\sigma}$ and $D_{\sigma{'}}$ for some $ \sigma,\sigma{'} \in \mathbb{Z}_3^2\setminus\left\lbrace \left( 0,0\right)  \right\rbrace $ and:
	\begin{enumerate}
		\item $ D_{\sigma}$ and $D_{\sigma{'}}$ are smooth at $ y $ and they meet transversally at $ y $;
		\item the inertia groups of $ D_{\sigma}$ and $D_{\sigma{'}}$ are different.
	\end{enumerate}

\end{enumerate}
\end{Proposition}

\noindent
Also from \cite[\rm Lemma 4.2, Proposition 4.2]{MR1103912} we have:
\begin{Proposition}\label{invariants of bitriple cover}
Let $ \xymatrix{f: X \ar[r]& Y} $ be a smooth $  \mathbb{Z}_{3}^2 $-cover with the building data $ \left\lbrace L_{\chi}, D_{\sigma} \right\rbrace_{\chi,\sigma}$. The surface $ X $ satisfies the following:
$$
\begin{adjustbox}{max width=\textwidth}
\begin{tabular}{l l l}
$ 3K_X $ &$ \equiv $ &$ f^*\left( 3K_Y + \sum\limits_{\sigma \ne \left( 0,0\right) } {2D_{\sigma} } \right) $	\\
$ K^2_X $ &$ = $ &$  \left( 3K_Y + \sum\limits_{\sigma \ne \left( 0,0\right) } {2D_{\sigma} } \right)^2 $ \\
$ p_g\left( X \right) $ &$ = $ &$ p_g\left( Y \right) +\sum\limits_{\chi  \ne \chi_{0,0} } {h^{0}\left( L_{\chi } + K_Y \right)}  $ \\
$ \chi\left( \mathcal{O}_X \right) $ &$ = $&$9\chi\left( \mathcal{O}_Y \right) + \sum\limits_{\chi \ne \chi_{0,0} } {\frac{1}{2}\left( L_{\chi } + K_Y \right)L_{\chi }}. $ 
\end{tabular}
\end{adjustbox}
$$
\end{Proposition}

\subsection{$ \mathbb{Z}_{3} $-coverings}\label{Z3-covering}
Let $ Y $ be a smooth surface with no non-trivial $ 3 $-torsion and let $ \xymatrix{f: X \ar[r] & Y} $ be a $ \mathbb{Z}_3 $-cover. The building data of this $ \mathbb{Z}_3 $-cover is $ \left\lbrace D_1, D_2, L_1, L_2\right\rbrace  $, where $ D_1, D_2, L_1, L_2 $ are divisors of $ Y $ such that $ D_1, D_2 $ are effective, $ L_1, L_2 $ are non-trivial and such that
\begin{align*}
3L_1 &\equiv D_1 + 2D_2\\
3L_2 &\equiv 2D_1 + D_2.
\end{align*}
\noindent
When the total branch locus $ B = D_1 + D_2 $  is smooth, the surface $ X $ is smooth and satisfies the following:
\begin{align*}
3K_X &\equiv f^{*}\left( 3K_Y + 2D_1 + 2D_2\right),\\
p_g\left( X\right) &= p_g\left( Y\right) + h^{0}\left( Y\right) + h^{0}\left( K_Y +L_1\right) + h^{0}\left( K_Y +L_2\right),\\
\chi\left( \mathcal{O}_X\right)  &= 3\chi\left( \mathcal{O}_Y\right) + \frac{1}{2}L_1\left( K_Y +L_1\right) +\frac{1}{2}L_1\left( K_Y +L_1\right).
\end{align*}

\noindent
In the case where the total branch locus $ B = D_1 + D_2 $ has only normal crossing singularities, from \cite[\rm Proposition 3.3]{MR1103912} the surface $ X $ has only rational singularities that are determined as follows:
\begin{Remark}\label{the type of singularities in the quotient surfaces}
Let $ Y $ be a smooth surface and let $ \xymatrix{f: X \ar[r] & Y} $ be a $ \mathbb{Z}_3 $-cover with the building data $ \left\lbrace D_1, D_2, L_1, L_2\right\rbrace  $. Assume that the total branch locus $ B = D_1 + D_2 $ has smooth components and only normal crossing singularities. Then
\begin{enumerate}
	\item if $ y $ is a double point in one of the $ D_i $, $ f^{-1}\left( y\right)  $ is a singular point of type $ A_2 $ in $ X $,
	\item if $ y $ belongs to $ D_1 \cap D_2$, $ f^{-1}\left( y\right)  $ is a singular point of type $ \frac{1}{3}\left( 1,1\right)  $ in $ X $.
\end{enumerate}	
\end{Remark}
\noindent
In Remark \ref{the type of singularities in the quotient surfaces}, a singular point of type $ A_2 $ is the quotient singularity whose minimal resolution is the union of two $ \left(-2\right) $-curves $ C_1, C_2$ such that $ C_1C_2=1 $ and a singular point of type $ \frac{1}{3}\left( 1,1\right)  $ is the quotient singularity whose minimal resolution is a smooth rational curve of self-intersection $ -3 $.\\

Let $ X $ be a smooth projective surface. Assume that the cyclic group $ \mathbb{Z}_3 $ acts on $ X $ with only isolated fixed points. By Cartan's lemma, the character of the representation of $ \mathbb{Z}_3 $ on the tangent space at an isolated fixed point is the sum of two non-trivial characters. The image of this point in the quotient surface $ X/{\mathbb{Z}_3} $ is a canonical singularity of type $ A_2 $ if the two characters are different and it is a singular point of type $ \frac{1}{3}\left( 1,1\right)  $ if the two characters are equal (see \cite[Proposition 5.3, Theorem 5.4]{MR2030225} and \cite[Lemma 2]{MR0084174}). 	    
\begin{Proposition}[{\cite[Proposition 2.1]{MR1960801}}]\label{singularities on Z_3 quotient}
Let $ X $ be a smooth projective surface. Assume that the group $ \mathbb{Z}_3 $ acts on $ X $ with only isolated fixed points. Denote by $ \tilde{X}_1 $ the minimal resolution of the quotient surface $ X/{\mathbb{Z}_3} $, then
\begin{align*}
K_{X}^2 &= 3K_{\tilde{X}_1}^2 +m\\
\chi\left( \mathcal{O}_{X} \right) &= 3\chi\left( \mathcal{O}_{\tilde{X}_1} \right) - \frac{2n+m}{3},
\end{align*}	    	    	
where $ m $ is the number of fixed points whose image in the quotient surface $ X/{\mathbb{Z}_3} $ is a singular point of type $ \frac{1}{3}\left( 1,1\right)  $ and $ n $ is the number of fixed points whose image in the quotient surface $ X/{\mathbb{Z}_3} $ is a singular point of type $ A_2 $.
\end{Proposition}    

Let $ P_1, P_2, \ldots, P_m $ be the fixed points for the $ \mathbb{Z}_3 $-action whose images in the quotient surface $ X/{\mathbb{Z}_3} $ are singularities of type $ \frac{1}{3}\left( 1,1\right)  $. We consider the following commutative diagram
$$
\xymatrix{ \overline{X} \ar[r]^{\overline{\pi}} \ar[d]^{\overline{f}}& X \ar[d]^{f}\\
\overline{X}_1 \ar[r]^{\overline{r}}& X/{\mathbb{Z}_3}
}
$$
where $ \overline{\pi} $ is the blow-up of $ X $ at the points $ P_1, P_2, \ldots, P_m $, the map $ \overline{r} $ is a partial resolution of $ X/{\mathbb{Z}_3} $ resolving only the singularities of type $ \frac{1}{3}\left( 1,1\right)  $, $ f $ is the quotient map with respect to the action of $ \mathbb{Z}_3 $, and $ \overline{f} $ is the quotient map with respect to the extended action $ \mathbb{Z}_3 $ on $ \overline{X} $. Denoting by $ E_1, E_2, \ldots, E_m $ the exceptional divisors corresponding to $ P_1, P_2, \ldots, P_m $, respectively, the canonical divisor of $ \overline{X} $ is
\begin{align*}
K_{\overline{X}} &\equiv \overline{\pi}^{*}\left( K_X\right) + \sum_{i=1}^{m}{E_i}.
\end{align*}
On the other hand, because $ \overline{f} $ is totally ramified on the $ \left( -3\right)$-curves 
$$ \overline{r}^{-1}\left( f\left( P_1\right) \right), \overline{r}^{-1}\left( f\left( P_2\right) \right), \ldots, \overline{r}^{-1}\left( f\left( P_m\right) \right),$$
by the Hurwitz formula, we obtain that
\begin{align*}
K_{\overline{X}} &\equiv \overline{f}^{*}\left( K_{\overline{X}_1}\right) + \sum_{i=1}^{m}{2E_i}. 
\end{align*}

\noindent
So we get that 
\begin{align*}
\overline{\pi}^{*}\left( K_X\right) \equiv \overline{f}^{*}\left( K_{\overline{X}_1}\right) + \sum_{i=1}^{m}{E_i}.
\end{align*}

Suppose that $ p_g\left( X\right) = p_g\left(\tilde{X}_1 \right)  $. Since the surface $ \overline{X}_1 $ has only canonical singularities $ A_2 $, we have that
\begin{align*}
h^{0}\left( \overline{\pi}^{*}\left( K_{X}\right) \right) = h^{0}\left( \overline{f}^{*}\left( K_{\overline{X}_1}\right) \right)
\end{align*}
Thus, the $ \left(-1 \right)  $-curves $ E_1, E_2, \ldots, E_m $ are fixed components of $ \overline{\pi}^{*}\left( K_{X}\right) $. Therefore the points $ P_1, P_2, \ldots, P_m $ are base points of $ \left|K_{X} \right|  $. We have thus proved the following proposition: 
\begin{Proposition}\label{base-point-of-triple-cover}
Let $ X $ be a smooth projective surface. Assume that the group $ \mathbb{Z}_3 $ acts on $ X $ with only isolated fixed points. Let $ P_1, P_2, \ldots, P_m $ be the fixed points for the $ \mathbb{Z}_3 $-action whose images in the quotient surface $ X/{\mathbb{Z}_3} $ are singular points of type $ \frac{1}{3}\left( 1,1\right)  $. Denote by $ \tilde{X}_1 $ the minimal resolution of the quotient surface $ X/{\mathbb{Z}_3} $. If 
\begin{align*}
p_g\left( X\right) = p_g\left(\tilde{X}_1 \right),
\end{align*}
the points $ P_1, P_2, \ldots, P_m $ are base points of $ \left|K_{X} \right|  $.
\end{Proposition}

\section{Constructions of the surfaces listed in Theorem \ref{THE MAIN THEOREM WITH D=3 PART 1}}
In this section, we construct the surfaces described in Theorem \ref{THE MAIN THEOREM WITH D=3 PART 1}. First, we recall the notion $ \mathbb{Z}_3 $-set of singularities which was introduced by S. Tan \cite{MR1141782}:
\begin{Definition}\label{Z_3-set}
	Let $ N $ be a set of finite singular points of a surface $ Y $ and let $ \xymatrix{\sigma : \tilde{Y} \ar[r] &Y} $ be the minimal resolution of $ Y $. We denote by $ C_1, C_2, \ldots, C_r $ the components of $ \sigma^{-1}\left( N\right)  $. If there exist integer numbers $ a_1, a_2, \ldots, a_r \in \left\lbrace 1,2\right\rbrace  $ and a divisor $ L $ of $ \tilde{Y} $ such that $ \sum_{i=1}^{r}{a_iC_i} \equiv 3L $, the set $ N $ is called a $ \mathbb{Z}_3 $-set of singularities on $ Y $.
\end{Definition}

\subsection{The main construction}\label{The main construction 1 of surfaces with d = 3}
The main idea of the construction of the surface $ X $ described in the first row of Theorem \ref{THE MAIN THEOREM WITH D=3 PART 1} is the following: we construct $ X $ as a $ \mathbb{Z}_3^2 $-cover $ \xymatrix{X \ar[r] & Y_3} $ of a del Pezzo surface $ Y_3 $ of degree $ 6 $ ramified on a total branch divisor $ B \in \left| 2f_1 + 3f_2 + 3f_3 \right| $ (See notations in Section \ref{Construction and computation of invariants part 1} ). We take the total branch locus $ B $ such that $ \chi\left( \mathcal{O}_{X} \right) = 6$, $ q\left( X\right) =0 $, that one of the $ \mathbb{Z}_3 $-quotients $ X_1 $ of $ X $ is a surface of general type with a $ \mathbb{Z}_3 $-set of $ 15 $ singular points of type $ A_2 $ and $ 6 $ singular points of type $ \frac{1}{3}\left( 1,1\right)  $, and that the following diagram commutes:
$$
\xymatrix{X \ar[0,3]^{\mathbb{Z}_3^2} \ar[1,1]^{\mathbb{Z}_3} \ar@{.>}[3,1]_{\varphi_{\left| K_X \right| }}&&& Y_3 \\
	&X_1 \ar[-1,2]^{\mathbb{Z}_3} &&\\
	&\overline{X}_1 \ar@{.>}[1,0]^{\varphi_{\left| K_{\overline{X}_1} \right| }} \ar[-1,0]^{\overline{r}}&&\\
	&\hskip 0.5cm \Sigma \subseteq \mathbb{P}^{4}&&} 	
$$

\subsubsection{Construction and computation of invariants}\label{Construction and computation of invariants part 1}
We denote by $ Y_3 $ the blow-up of $ \mathbb{P}^2$ at three distinct non-collinear points $ P_1 $, $ P_2 $, $ P_3 $. Let us denote by 
\begin{center}
	\begin{tabular}{l l}
		$ l $: & the pullback of a general line in $ \mathbb{P}^2$,\\
		$ e_1 $, $ e_2 $, $ e_3 $:&the exceptional divisors corresponding to $ P_1 $, $ P_2 $, $ P_3 $, respectively,\\
		$ f_1 $, $ f_2$, $ f_3$:&the strict transforms of a general line through $ P_1 $, $ P_2 $, $ P_3 $, respectively,\\
		$ h_{12} $, $ h_{23} $, $ h_{31} $:&the strict transforms of the lines $ P_1  P_2 $, $ P_2  P_3 $, $ P_3  P_1 $, respectively.
	\end{tabular}	
\end{center}
The anti-canonical class $ -K_{Y_3} \equiv f_1 + f_2 + f_3 $ is very ample and the linear system $ \left| -K_{Y_3} \right|  $ embeds $ Y_3 $ as a smooth del Pezzo surface of degree $ 6 $ in $ \mathbb{P}^6 $.\\		  

We consider the following smooth divisors:
\begin{align*}
D_{01}:= f_{11}, D_{02}:= f_{12}, D_{22}:= f_{21}+f_{22}+f_{23}, D_{12}:= f_{31}+f_{32}+f_{33}
\end{align*}
\noindent
and $ D_{ij} = 0 $ for the other $ \left( i,j\right)  $, where $ f_{11}, f_{12} \in \left| f_1\right| $, $ f_{21}, f_{22}, f_{23} \in \left| f_2\right| $, $ f_{31}, f_{32}, f_{33} \in \left| f_3\right| $ are distinct divisors of $ Y_3 $ such that no more than two of these divisors go through the same point. We consider the following non-trivial divisors:
$$	
\begin{tabular}{r r r r}
$ L_{10} :=$ & $$&$f_2 $&$+2f_3  $ \\
$ L_{01} :=$ & $f_1$&$+f_2 $&$+f_3 $\\
$ L_{20} :=$ &  $$&$2f_2 $&$+f_3 $\\
$ L_{02} :=$ & $f_1$&$+2f_2 $&$+2f_3 $\\		  
$ L_{11} :=$ &  $f_1$&$+2f_2 $&$ $\\
$ L_{22} :=$ &  $f_1$&$+f_2 $&$ $\\
$ L_{12} :=$ &  $f_1$&$ $&$+f_3 $\\
$ L_{21} :=$ &  $f_1$&$ $&$+2f_3 $.
\end{tabular} 
$$
\noindent
These divisors $ D_{\sigma} $, $ L_{\chi} $ satisfy the following relations:
$$
\begin{adjustbox}{max width=\textwidth}
\begin{tabular}{l l r l r l r l r r r r r}
$ 3L_{10} $&$ \equiv $&&$   $& &$ $&$D_{22} $&$ +$&$2D_{12} $ &$ \equiv $ & $$&$3f_2 $&$+6f_3  $\\
$ 3L_{01} $&$ \equiv $ & $D_{01}$&$ +$ & $2D_{02}  $&$ +$&$D_{22} $&$ +$&$D_{12} $ &$ \equiv $ & $3f_1$&$+3f_2 $&$+3f_3 $\\
$ 3L_{20} $&$ \equiv $&&$ $&&$ $&$2D_{22} $&$ +$&$D_{12} $ &$ \equiv $ &  $$&$6f_2 $&$+3f_3 $\\
$ 3L_{02} $&$ \equiv$ & $2D_{01} $&$ +$ & $D_{02} $&$ +$&$2D_{22} $&$ +$&$2D_{12} $ &$ \equiv $ & $3f_1$&$+6f_2 $&$+6f_3 $\\
$ 3L_{11} $&$ \equiv $ & $D_{01}$&$ +$ & $2D_{02}  $  &$ +$&$2D_{22} $&& &$ \equiv $ &  $3f_1$&$+6f_2 $&$ $\\
$ 3L_{22} $&$ \equiv $ & $2D_{01}$&$ +$ & $D_{02}  $  &$ +$&$D_{22} $&& &$ \equiv $ &  $3f_1$&$+3f_2 $&$ $\\
$ 3L_{12} $&$ \equiv $ & $2D_{01}$&$ +$ & $D_{02}  $&$  $&&$ +$&$D_{12} $ &$ \equiv $ &  $3f_1$&$ $&$+3f_3 $\\
$ 3L_{21} $&$ \equiv $ & $D_{01}$&$ +$ & $2D_{02}  $&$  $&&$ +$&$2D_{12} $&$ \equiv $ &  $3f_1$&$ $&$+6f_3 $.
\end{tabular}
\end{adjustbox}
$$
\noindent
Thus by Proposition \ref{Construction of bitriple cover}, the divisors $ D_{\sigma} $, $ L_{\chi} $ define a $ \mathbb{Z}^2_3 $-cover $ \xymatrix{g: X \ar[r] & Y_3}  $. Moreover, by Propositions \ref{Smoothness condition of bitriple cover} and \ref{invariants of bitriple cover}, the surface $ X $ is smooth and satisfies the following:
\begin{align*}
3K_X &\equiv g^*\left( f_1 + 3f_2 + 3f_3\right). 
\end{align*}		  
We notice that a surface is of general type and minimal if the canonical divisor is big and nef (see e.g. \cite[\rm Section 2]{MR2931875}). Because the divisor $ 3K_X $ is the pull-back of a big and nef divisor, the canonical divisor $ K_X $ is big and nef. Thus, the surface $ X $ is minimal and of general type. Furthermore, by Proposition \ref{invariants of bitriple cover}, the surface $ X $ possesses the following invariants:
\begin{align*}
K_X^2& =\left( f_1 + 3f_2 + 3f_3\right)^2= 30,\\
p_g\left( X\right) &=5, \chi\left( \mathcal{O}_X\right) =6, q\left( X\right)  = 0.   
\end{align*}

\subsubsection{The canonical map and the canonical image}
Let $ X $ be as above and let $ \Sigma $ denote the image of the canonical map $ \varphi_{\left| K_X \right| } $. We show that the degree $d$ of the canonical  map is $ 3 $ and that the canonical image $ \Sigma $ is a surface of general type. We consider the cyclic subgroup $ \Gamma:= \left\langle \left( 1,0\right) \right\rangle  $ of $ \mathbb{Z}_3^2 $. Because $ h^0\left( L_{ \chi} + K_{Y_3} \right)  = 0 $ for all $ \chi \notin \Gamma^\perp $, by Proposition \ref{Factorization of the canonical map through a quotient surface}, the canonical map $ \varphi_{\left| K_X \right| } $ is essentially the composition of the quotient map $g_{2}: \xymatrix{X \ar[r]& X_{1}:= X/\left\langle \left( 1,0\right) \right\rangle } $ with the canonical map $\varphi_{\left| K_{\tilde{X}_{1}} \right| }$ of the minimal resolution $ \tilde{X}_{1} $ of $ X_1 $. 

The surface $ X_1 $ is the $ \mathbb{Z}_3 $-cover $ \xymatrix{g_1 : X_1 \ar[r]& Y_3} $ with the following building data (see Section \ref{Z3-covering}):
\begin{align*}
D_{1}&:= D_{01} + D_{22} +D_{12}, &D_{2}&: = D_{02}\\		  
L_{1}&: = L_{01}, &L_{2}&: = L_{02}.		  
\end{align*}
\noindent
Because the divisor $ D_1 $ has fifteen double points and $ D_1D_2 = 6 $, by Remark \ref{the type of singularities in the quotient surfaces} the quotient surface $ X_1 $ has fifteen singular points of type $ A_2 $ and six singular points of type $ \frac{1}{3}\left( 1,1\right)  $. We consider the following maps
$$
\xymatrix{\overline{X}_1 \ar[r]^{\overline{g}_1} & \overline{Y}_3 \ar[r]^{\pi_1}& Y_3}
$$
where the map $ \pi_{1} $ is the blow-up at the intersection points between $ D_{1} $ and $ D_{2} $ (i.e. between $ D_{02} $ and $ D_{22} +D_{12}$), and the map $ \overline{g}_1 $ is the $ \mathbb{Z}_3 $-cover which is ramified on the strict transforms of $ D_1, D_2 $. The surface $ \overline{X}_1 $ is a partial resolution of $ X_1 $ where the singularities of type $ \frac{1}{3}\left( 1,1\right)  $ are resolved. The surface $ \overline{X}_1 $ has fifteen singular points of type $ A_2 $ and satisfies the following
\begin{align*}
3K_{\overline{X}_1} &\equiv \overline{g}_1^{*}\left( \pi_1^{*}\left( f_1 + 3f_2 + 3f_3\right)  -\sum_{i=1}^{6}{h_i}\right), 
\end{align*}
where $ h_1, h_2, \ldots, h_6 $ are the exceptional divisors of $ \pi_1 $. Let $ \tilde{X}_1 $ be the minimal resolution of $ X_1 $. Because the surface $ \overline{X}_1 $ has only canonical singularities,  $ K_{\tilde{X}_1}^2 = K_{\overline{X}_1}^2 = \left( \pi_1^{*}\left( f_1 + 3f_2 + 3f_3\right)  -\sum_{i=1}^{6}{h_i}\right)^2 =8 $ and
\begin{align*}
p_g\left( \tilde{X}_1\right) =5, \chi\left( \mathcal{O}_{\tilde{X}_1}\right) =6, q\left( \tilde{X}_1\right)  = 0. 
\end{align*}		  
\noindent
Moreover, the surface $ \tilde{X}_1 $ satisfies the following
\begin{align*}
3K_{\tilde{X}_1} &\equiv \tilde{g}_1^{*}\left( \pi_1^{*}\left( f_1 + 3f_2 + 3f_3\right) -\sum_{i=1}^{6}{h_i}\right), 
\end{align*}
\noindent
where $ \tilde{g}_1 $ is the composition of the map $ \overline{g}_1 $ and the resolution map which resolves the singular points of type $ A_2 $.

Remark that the surface $ \tilde{X}_1 $ has no non-trivial $ 3$-torsion. In fact, assume that $ \tilde{X}_1 $ has a non-trivial divisor $ T $ such that $ 3T\equiv 0 $. We consider the \'{e}tale cover $ \xymatrix{h: Z \ar[r]& \tilde{X}_1} $ of degree $ 3 $ which is ramified on $ 3T $. So $ Z $ is a minimal surface with
\begin{align*}
K_Z^2 = 24, \chi\left( \mathcal{O}_Z \right) = 18. 
\end{align*}
This contradicts Noether's inequality. So we get
\begin{align*}
K_{\tilde{X}_1} &\equiv \theta +\tilde{g}_1^{*}\left( \pi_1^{*}\left( f_2\right) \right) + \tilde{g}_1^{*}\left( \pi_1^{*}\left( f_3\right) \right)
\end{align*}
where $ \theta:=\tilde{g}_1^{*}\left( \pi^{*}_1\left( f_{12}\right)  -\sum\limits_{i=1}^{6}{h_i}\right)_{\text{red}} $. We remark that the divisor $ \pi_1^{*}\left(f_{12}\right)  -\sum\limits_{i=1}^{6}{h_i} $ is effective and is contained in the branch locus. So the divisor $ \theta $ is a $ \left( -2\right)  $-curve.

Notice that a general curve in $ \left| \tilde{f}_2 \right| $ is not hyper-elliptic, where $ \tilde{f}_2:= \tilde{g}_1^{*}\left( \pi^{*}_1\left( f_2\right) \right)$. Indeed, a general curve in $ \left| \tilde{f}_2 \right| $ has a base point free $ g_3^1 $ which is a special linear system. So a general curve in $ \left| \tilde{f}_2 \right| $ is not hyper-elliptic because a special linear system on a hyper-elliptic curve is always composed of a multiple of the $ g_2^1 $ plus some base points (see \cite[ \rm Lemma 5]{MR875340}).

Moreover, the restriction map $ \xymatrix{H^{0}\left( \tilde{X}_1,\mathcal{O}_{\tilde{X}_1} \left( K_{\tilde{X}_1} \right) \right) \ar[r] & H^{0}\left( \tilde{f}_2, \mathcal{O}_{\tilde{f}_2}\left( K_{\tilde{X}_1}\right)\right)  }  $ is surjective. In fact, we consider the following long exact sequences
$$
\begin{adjustbox}{max width=\textwidth}
$$
\xymatrix{0 \ar[r]& H^{0}\left( \tilde{X}_1,\mathcal{O}_{\tilde{X}_1} \left( K_{\tilde{X}_1} -\tilde{f}_2 \right) \right) \ar[r] &  H^{0}\left( \tilde{X}_1,\mathcal{O}_{\tilde{X}_1} \left( K_{\tilde{X}_1} \right) \right) \ar[r] & H^{0}\left( \tilde{f}_2, \mathcal{O}_{\tilde{f}_2}\left( K_{\tilde{X}_1}\right)\right) \ar[r]& \ldots }
$$ 
\end{adjustbox}
$$
\noindent
Because $ h^{0}\left( \tilde{X}_1,\mathcal{O}_{\tilde{X}_1} \left( K_{\tilde{X}_1} \right) \right)=5 $, $ h^{0}\left( \tilde{f}_2, \mathcal{O}_{\tilde{f}_2}\left( K_{\tilde{X}_1}\right)\right)=3 $ and
\begin{align*}
h^{0}\left( \tilde{X}_1,\mathcal{O}_{\tilde{X}_1} \left( K_{\tilde{X}_1} -\tilde{f}_2 \right) \right)&=h^{0}\left( \tilde{X}_1,\mathcal{O}_{\tilde{X}_1} \left( \theta +\tilde{f}_3 \right) \right)= h^{0}\left( \tilde{X}_1,\mathcal{O}_{\tilde{X}_1} \left( \tilde{f}_3 \right) \right) =2,
\end{align*}
\noindent
the restriction map $ \xymatrix{H^{0}\left( \tilde{X}_1,\mathcal{O}_{\tilde{X}_1} \left( K_{\tilde{X}_1} \right) \right) \ar[r] & H^{0}\left( \tilde{f}_2, \mathcal{O}_{\tilde{f}_2}\left( K_{\tilde{X}_1}\right)\right)  }  $ is surjective.

On the other hand, the linear system $ \left| \tilde{f}_2 \right| $ moves in a rational pencil. So $ \varphi_{\left| K_{\tilde{X}_1} \right| } $ separates the fibers of $ \left| \tilde{f}_2 \right|  $. This yields that  $ \varphi_{\left| K_{\tilde{X}_1} \right| } $ is birational.\\

We now show that the linear system $ \left| K_{\tilde{X}_1} \right| $ is base point free. We notice that the divisor $ \theta $ is a $ \left( -2\right)  $-curve of $ \tilde{X}_1 $ and 
\begin{align*}
h^{0}\left( \tilde{X}_1, \tilde{g}_1^{*}\left( \pi_1^{*}\left( f_2\right) \right) + \tilde{g}_1^{*}\left( \pi_1^{*}\left( f_3\right) \right)  \right) =4.
\end{align*}
\noindent
Because $ h^{0}\left(\tilde{X}_1, K_{\tilde{X}_1} \right)  \ne h^{0}\left( \tilde{X}_1, \tilde{g}_1^{*}\left( \pi_1^{*}\left( f_2\right) \right) + \tilde{g}_1^{*}\left( \pi_1^{*}\left( f_3\right) \right)  \right) $, $ \theta K_{\tilde{X}_1}=0 $ and $ \left| \tilde{g}_1^{*}\left( \pi_1^{*}\left( f_2\right) \right) + \tilde{g}_1^{*}\left( \pi_1^{*}\left( f_3\right) \right) \right|  $ is base point free, it follows from
\begin{align*}
K_{\tilde{X}_1} \equiv \theta +\tilde{g}_1^{*}\left( \pi_1^{*}\left( f_2\right) \right) + \tilde{g}_1^{*}\left( \pi_1^{*}\left( f_3\right) \right)
\end{align*}
\noindent
that $ \left| K_{\tilde{X}_1} \right| $ is base point free. Therefore, the canonical map $ \varphi_{\left| K_{\tilde{X}_1} \right| }$ is a birational morphism and the image $ \Sigma $ is a complete intersection of a quadric and a quartic in $ \mathbb{P}^4 $ (see \cite[Theorem 4.6]{MR938031} and \cite[Theorem 1.5]{MR1076174}).

The second cover $ g_2 $ is ramified on the $ \mathbb{Z}_3 $-set of fifteen singular points of type $ A_2 $ and six singular points of type $ \frac{1}{3}\left( 1,1\right)  $. Because $ \varphi_{\left| K_{\tilde{X}_1} \right| } $ is birational, $ \varphi_{\left| K_{X} \right| } $ is a $ 3$-to-$1 $ map onto a surface of general type. Thus, we obtain the surface described in the first row of Theorem \ref{THE MAIN THEOREM WITH D=3 PART 1}.\\

Finally, because the $ \mathbb{Z}_3$-cover $ \xymatrix{X \ar[r]^{g_2} & X_1} $ is ramified on the $ \mathbb{Z}_3 $-set of six singular points of type $ \frac{1}{3}\left( 1,1\right)  $ and fifteen singular points of type $ A_2 $, and $ p_g\left( X\right)  = p_g\left( \tilde{X}_1\right) $,  by Proposition \ref{base-point-of-triple-cover}, the canonical system $ \left| K_X \right|  $ has six base points.

\subsection{Variations}
We now consider variations of the construction of Section \ref{The main construction 1 of surfaces with d = 3}, to obtain the other surfaces in Theorem \ref{THE MAIN THEOREM WITH D=3 PART 1}. In order to achieve this, we impose ordinary triple points on the branch locus of the $ \mathbb{Z}_3^2 $-cover in Section \ref{The main construction 1 of surfaces with d = 3} and then we resolve the singularities. This is similar to the work of S. L. Tan for triple covers \cite{MR1150664}. We require that $ D_{\sigma_1} $, $ D_{\sigma_2} $ and $ D_{\sigma_3} $ have different inertia groups and meet transversally at a point. There are two possibilities:		  
\begin{enumerate}
	\item The coefficients of $ D_{\sigma_1} $, $ D_{\sigma_2} $, $ D_{\sigma_3} $ are different in the relation of $ 3L_{01} $ (in Proposition \ref{Construction of bitriple cover});
	\item The coefficients of $ D_{\sigma_1} $, $ D_{\sigma_2} $, $ D_{\sigma_3} $ are the same in the relation of $ 3L_{01} $ (in Proposition \ref{Construction of bitriple cover}).
\end{enumerate}

\subsubsection{Variation 1}
The first variation of the construction in Section \ref{The main construction 1 of surfaces with d = 3} consists of allowing the three divisors $ D_{02} $, $ D_{22} $, $ D_{12} $ meet transversally at a point $ P $. We notice that the coefficients of $ D_{02} $, $ D_{22} $, $ D_{12} $ are different in the relation of $ 3L_{01} $ (in Proposition \ref{Construction of bitriple cover}). More precisely, we take the branch locus 
\begin{align*}
D_{01} \in \left| f_1\right|, D_{02} \in \left| f_1\right|, D_{22} \in \left| 3f_2\right|, D_{12} \in \left| 3f_3\right|
\end{align*}
\noindent
and $ D_{ij} = 0 $ for the other $ \left( i,j\right)  $ as Section \ref{The main construction 1 of surfaces with d = 3} except that now we require that the three divisors $ D_{02} $, $ D_{22} $, $ D_{12} $ meet transversally at a point $ P $. If we take a $ \mathbb{Z}_{3}^2 $-cover which is ramified on this branch locus, we obtain a surface with singularities. The singularities appearing in this case are resolved as follows: since the strict transforms of $ D_{01} $, $ D_{02} $, $ D_{22} $, $ D_{12} $ of the blow-up at $ P $ do not satisfy the condition of Proposition \ref{Construction of bitriple cover}, we can not take a $ \mathbb{Z}_{3}^2 $-cover which is ramified on these strict transforms. To solve this problem, we add the exceptional divisor to $ D_{01} $. This means the new branch curve $ D_{01} $ is the strict transform $ \overline{D}_{01} $ of $ D_{01} $ and the exceptional divisor $ E $. The new branch curves $ D_{\sigma} $ are the strict transforms $ \overline{D}_{\sigma} $ of $ D_{\sigma} $ for all $ D_{\sigma} \ne D_{01} $. These new branch curves now satisfy the condition of Proposition \ref{Construction of bitriple cover}. We, however, still get a surface with singularities because the new branch curves $ D_{01} $ and $ D_{02} $ have a common point. Notice that $ D_{01} $ and $ D_{02} $ have the same inertia group. To fulfill the condition of Proposition \ref{Smoothness condition of bitriple cover}, we blow up the intersection point between the exceptional divisor and the strict transform of $ D_{02} $. The preimage of $ P $ is the union of four rational curves $ C_0, C_{1}, C_{2}, C_{3} $ such that $ C_0^2 = -2 $, $ C_{i}^2 = -3 $, $ C_0C_{i}=1 $ for all $ i \in \left\lbrace 1,2,3\right\rbrace  $ and that $ C_{i}C_{j}=0 $ for all $ i \ne j $.\\

We now present this variation in details. Let us denote by $ Y_4 $ the blow-up of $ \mathbb{P}^2$ at four points in general position $ P_1, P_2, P_3, P_4 $ and denote by 
\begin{center}
	\begin{tabular}{l l}
		$ l $: & the pullback of a general line in $ \mathbb{P}^2$,\\
		$ e_1 $, $ e_2 $, $ e_3 $, $ e_4 $:&the exceptional divisors corresponding to $ P_1 $, $ P_2 $, $ P_3 $, $ P_4 $, respectively,\\
		$ f_1 $, $ f_2$, $ f_3$, $ f_4$:&the strict transforms of a general line through $ P_1 $, $ P_2 $, $ P_3 $, $ P_4 $, respectively,\\
		$ h_{ij} $:&the strict transforms of the line $ P_i  P_j $, for all $ i, j \in \left\lbrace 1,2,3,4\right\rbrace  $, respectively.
	\end{tabular}	
\end{center}
The anti-canonical class 
\begin{align*}
-K_{Y_4} &\equiv f_1 + f_2 + f_3 - e_4 \\
&\equiv f_1 + f_2 + f_4 - e_3 \\
&\equiv f_1 + f_3 + f_4 - e_2 \\
&\equiv f_2 + f_3 + f_4 - e_1  
\end{align*}
\noindent
is very ample and the linear system $ \left| -K_{Y_4} \right|  $ embeds $ Y_4 $ as a smooth del Pezzo surface of degree $ 5 $ in $ \mathbb{P}^5 $.\\

Let $ \overline{Y}_4 $ be the blow-up of $ Y_4 $ at the intersection point between $ e_4 $ and $ h_{14} $. We denote by
\begin{center}
	\begin{tabular}{l l}
		$ \overline{e}'_4 $: & the exceptional divisor,\\
		$ \overline{e}_4 $, $ \overline{h}_{14} $:&the strict transforms of $ e_4 $, $ h_{14} $, respectively.\\
		$ \overline{e}_i $, $ \overline{f}_i $, $ \overline{h}_{ij} $, $ \overline{l} $:&can also be viewed as strict transforms of $ e_i $, $ f_i $, $ h_{ij} $, $ l $, respectively except $ e_4 $ and $ h_{14} $.
	\end{tabular}	
\end{center}

We consider the following smooth divisors:
\begin{align*}
D_{01}:= \overline{f}_{11} + \overline{e}_4, D_{02}:= \overline{h}_{14}, D_{22}:= \overline{h}_{24}+\overline{f}_{21}+\overline{f}_{22}, D_{12}:= \overline{h}_{34}+\overline{f}_{31}+\overline{f}_{32}
\end{align*}
\noindent
and $ D_{ij} = 0 $ for the other $ \left( i,j\right)  $, where $ \overline{f}_{11} \in \left| \overline{f}_1\right| $, $ \overline{f}_{21}, \overline{f}_{22} \in \left| \overline{f}_2\right| $, $ \overline{f}_{31}, \overline{f}_{32} \in \left| \overline{f}_3\right| $ are distinct divisors of $ Y_4 $ such that no more than two of these divisors go through the same point. We remark that 
\begin{align*}
\overline{h}_{ij} &\equiv \overline{f}_i - \overline{e}_i - \overline{e}_j \hskip 1cm \text{ for all } \hskip 0.5cm \left( i,j\right) \ne \left( 1,4\right),\\
\overline{h}_{14} &\equiv \overline{f}_1 - \overline{e}_1 - \overline{e}_4 -\overline{e}'_4.
\end{align*} 

\noindent
We have that
\begin{align*}
D_{01} &\equiv \overline{f}_1 + \overline{e}_4, &D_{02} &\equiv \overline{f}_1-\overline{e}_1 - \overline{e}_4 -2\overline{e}'_4\\
D_{22} &\equiv 3\overline{f}_2 - \overline{e}_4 -\overline{e}'_4, &D_{12} &\equiv 3\overline{f}_3 - \overline{e}_4 -\overline{e}'_4.
\end{align*}

\noindent
In addition, we consider the following non-trivial divisors:
$$	
\begin{tabular}{r r r r r r}
$ L_{10} :=$ & $$&$\overline{f}_2 $&$+2\overline{f}_3  $&$ -\overline{e}_4 $& $ -\overline{e}'_4 $\\
$ L_{01} :=$ & $\overline{f}_1$&$+\overline{f}_2 $&$+\overline{f}_3 $&$ -\overline{e}_4 $&$ -2\overline{e}'_4 $\\
$ L_{20} :=$ &  $$&$2\overline{f}_2 $&$+\overline{f}_3 $&$ -\overline{e}_4 $&$ -\overline{e}'_4 $\\
$ L_{02} :=$ & $\overline{f}_1$&$+2\overline{f}_2 $&$+2\overline{f}_3 $&$ -\overline{e}_4 $&$ -2\overline{e}'_4 $\\		  
$ L_{11} :=$ &  $\overline{f}_1$&$+2\overline{f}_2 $&$ $&$ -\overline{e}_4 $&$ -2\overline{e}'_4 $\\
$ L_{22} :=$ &  $\overline{f}_1$&$+\overline{f}_2 $&$ $&&$ -\overline{e}'_4 $\\
$ L_{12} :=$ &  $\overline{f}_1$&$ $&$+\overline{f}_3 $&&$ -\overline{e}'_4 $\\
$ L_{21} :=$ &  $\overline{f}_1$&$ $&$+2\overline{f}_3 $&$ -\overline{e}_4 $&$ -2\overline{e}'_4 $.
\end{tabular} 
$$
\noindent
The divisors $ D_{\sigma} $, $ L_{\chi} $ satisfy the following relations:
$$
\begin{adjustbox}{max width=\textwidth}
\begin{tabular}{l l r l r l r l r r r r r r r}
$ 3L_{10} $&$ \equiv $&&$   $& &$ $&$D_{22} $&$ +$&$2D_{12} $&$ \equiv$ & $$&$3\overline{f}_2 $&$+6\overline{f}_3  $&$ -3\overline{e}_4 $& $ -3\overline{e}'_4 $\\
$ 3L_{01} $&$ \equiv $ & $D_{01}$&$ +$ & $2D_{02}  $&$ +$&$D_{22} $&$ +$&$D_{12} $&$ \equiv$ & $3\overline{f}_1$&$+3\overline{f}_2 $&$+3\overline{f}_3 $&$ -3\overline{e}_4 $&$ -6\overline{e}'_4 $\\
$ 3L_{20} $&$ \equiv $&&$ $&&$ $&$2D_{22} $&$ +$&$D_{12} $ &$ \equiv$ &  $$&$6\overline{f}_2 $&$+3\overline{f}_3 $&$ -3\overline{e}_4 $&$ -3\overline{e}'_4 $\\
$ 3L_{02} $&$ \equiv$ & $2D_{01} $&$ +$ & $D_{02} $&$ +$&$2D_{22} $&$ +$&$2D_{12} $&$ \equiv$ & $3\overline{f}_1$&$+6\overline{f}_2 $&$+6\overline{f}_3 $&$ -3\overline{e}_4 $&$ -6\overline{e}'_4 $\\
$ 3L_{11} $&$ \equiv $ & $D_{01}$&$ +$ & $2D_{02}  $  &$ +$&$2D_{22} $&&&$ \equiv$ &  $3\overline{f}_1$&$+6\overline{f}_2 $&$ $&$ -3\overline{e}_4 $&$ -6\overline{e}'_4 $\\
$ 3L_{22} $&$ \equiv $ & $2D_{01}$&$ +$ & $D_{02}  $  &$ +$&$D_{22} $&&&$ \equiv$ &  $3\overline{f}_1$&$+3\overline{f}_2 $&$ $&&$ -3\overline{e}'_4 $\\
$ 3L_{12} $&$ \equiv $ & $2D_{01}$&$ +$ & $D_{02}  $&$  $&&$ +$&$D_{12} $&$ \equiv$ &  $3\overline{f}_1$&$ $&$+3\overline{f}_3 $&&$ -3\overline{e}'_4 $\\
$ 3L_{21} $&$ \equiv $ & $D_{01}$&$ +$ & $2D_{02}  $&$  $&&$ +$&$2D_{12} $ &$ \equiv$ &  $3\overline{f}_1$&$ $&$+6\overline{f}_3 $&$ -3\overline{e}_4 $&$ -6\overline{e}'_4 $.
\end{tabular}
\end{adjustbox}
$$
\noindent
Thus, by Proposition \ref{Construction of bitriple cover}, the divisors $ D_{\sigma} $, $ L_{\chi} $ define a $ \mathbb{Z}^2_3 $-cover $ \xymatrix{g: X \ar[r] & \overline{Y}_4}  $. Furthermore, by Propositions \ref{Smoothness condition of bitriple cover} and \ref{invariants of bitriple cover}, the surface $ X $ is smooth and satisfies the following:
\begin{align*}
3K_X &\equiv g^*\left( \overline{f}_1 + 3\overline{f}_2 + 3\overline{f}_3 -\overline{e}_4 -2\overline{e}'_4\right). 
\end{align*}		  
Since $ 3K_X $ is the pull-back of a nef and big divisor, the canonical divisor $ K_X $ is nef and big. Thus, the surface $ X $ is minimal and of general type. Furthermore, by Proposition \ref{invariants of bitriple cover}, the surface $ X $ possesses the following invariants:
\begin{align*}
K_X^2&=\left( \overline{f}_1 + 3\overline{f}_2 + 3\overline{f}_3 -\overline{e}_4 -2\overline{e}'_4\right)^2= 28,\\
p_g\left( X\right) &=5, \chi\left( \mathcal{O}_X\right) =6, q\left( X\right)  = 0.   
\end{align*}

As in Section \ref{The main construction 1 of surfaces with d = 3}, because the divisor $ D_{01} + D_{22} + D_{12} $ has sixteen double points and $ D_{02}\left( D_{01} + D_{22} + D_{12}\right) = 4 $, by Remark \ref{the type of singularities in the quotient surfaces} the quotient surface $ X_{1}:= X/\left\langle \left( 1,0\right) \right\rangle$ has sixteen singular points of type $ A_2 $ and four singular points of type $ \frac{1}{3}\left( 1,1\right)  $. We consider the following maps
$$
\xymatrix{\overline{X}_1 \ar[r]^{\overline{g}_1} & \hat{Y}_4  \ar[r]^{\pi_1}& \overline{Y}_4}
$$
where the map $ \pi_{1} $ is the blow-up at the intersection points between $ D_{02} $ and $ D_{22} +D_{12}$, the map $ \overline{g}_1 $ is the $ \mathbb{Z}_3 $-cover which is ramified on the strict transforms of $ D_{02}, D_{01}+D_{22}+D_{12} $. The surface $ \overline{X}_1 $ is a partial resolution of $ X_1 $ where the singularities of type $ \frac{1}{3}\left( 1,1\right)  $ are resolved. The surface $ \overline{X}_1 $ has sixteen singular points of type $ A_2 $ and satisfies the following
\begin{align*}
3K_{\overline{X}_1} &\equiv \overline{g}_1^{*}\left( \pi_1^{*}\left( \overline{f}_1 + 3\overline{f}_2 + 3\overline{f}_3\right) -\overline{e}_4 -2\overline{e}'_4 -\sum_{i=1}^{4}{h_i}\right), 
\end{align*}
where $ h_1, h_2, h_3, h_4 $ are the exceptional divisors of $ \pi_1 $. Because the surface $ \overline{X}_1 $ has only canonical singularities, the minimal resolution $ \tilde{X}_1 $ of $ X_1 $ has the following invariants:
\begin{align*}
K_{\tilde{X}_1}^2 &= K_{\overline{X}_1}^2 =\left( \pi_1^{*}\left( \overline{f}_1 + 3\overline{f}_2 + 3\overline{f}_3\right) -\overline{e}_4 -2\overline{e}'_4 -\sum_{i=1}^{4}{h_i}\right)^2= 8,\\
p_g\left( \tilde{X}_1\right) &=5, \chi\left( \mathcal{O}_{\tilde{X}_1}\right) =6, q\left( \tilde{X}_1\right)  = 0. 
\end{align*}		  
\noindent
Moreover, the surface $ \tilde{X}_1 $ satisfies the following
\begin{align*}
3K_{\tilde{X}_1} &\equiv \tilde{g}_1^{*}\left( \pi_1^{*}\left( \overline{f}_1 + 3\overline{f}_2 + 3\overline{f}_3\right) -\overline{e}_4 -2\overline{e}'_4 -\sum_{i=1}^{4}{h_i}\right), 
\end{align*}
\noindent
where $ \tilde{g}_1 $ is the composition of the map $ \overline{g}_1 $ and the resolution map which resolves the singular points of type $ A_2 $. Since the surface $ \tilde{X}_1 $ has no non-trivial $ 3$-torsion, the canonical divisor of $ \tilde{X}_1 $ is
\begin{align*}
K_{\tilde{X}_1} &\equiv \theta +\tilde{g}_1^{*}\left( \pi_1^{*}\left( \overline{f}_2\right) \right) + \tilde{g}_1^{*}\left( \pi_1^{*}\left( \overline{f}_3\right) \right)
\end{align*}
where $ \theta:=\tilde{g}_1^{*}\left( \pi_1^{*}\left( h_{14}-\overline{e}_4 -2\overline{e}'_4\right)  -\sum\limits_{i=1}^{4}{h_i}\right)_{\text{red}} $. We remark that the divisor $ \pi_1^{*}\left( h_{14}-\overline{e}_4 -2\overline{e}'_4\right)  -\sum\limits_{i=1}^{4}{h_i} $ is effective and is contained in the branch locus. So the divisor $ \theta $ is a $ \left( -2\right)  $-curve.

Analogously to Section \ref{The main construction 1 of surfaces with d = 3}, the canonical map $ \varphi_{\left| K_{\tilde{X}_1} \right| } $ is birational. Thus the canonical map $ \varphi_{\left| K_{X} \right| } $ is a $ 3$-to-$1 $ map onto a surface of general type. Furthermore, because the $ \mathbb{Z}_3$-cover $ \xymatrix{X \ar[r]^{g_2} & X_1} $ is ramified on the $ \mathbb{Z}_3 $-set of four singularities of type $ \frac{1}{3}\left( 1,1\right)  $ and sixteen singularities of type $ A_2 $ and $ p_g\left( X\right)  = p_g\left( \tilde{X}_1\right) $,  by Proposition \ref{base-point-of-triple-cover}, the canonical system $ \left| K_X \right|  $ has four base points. So we obtain the surface described in the second row of Theorem \ref{THE MAIN THEOREM WITH D=3 PART 1}.

\begin{Remark}
	Similarly as above, if we ask the three branch components $ D_{02} $, $ D_{22} $, $ D_{12} $ meet transversally at two (three) distinct points, we obtain the surface described in the third (the fourth) row of Theorem \ref{THE MAIN THEOREM WITH D=3 PART 1}. 
\end{Remark}

\subsubsection{Variation 2}
In the second variation of the construction given in Section Section \ref{The main construction 1 of surfaces with d = 3}, we again impose an ordinary triple point to be the branch locus but in a different way.  We require that the three branch components $ D_{01} $, $ D_{22} $ and $ D_{12} $ meet transversally at a point $ P $. The coefficients of $ D_{01} $, $ D_{22} $ and $ D_{12} $ are the same in the relation of $ 3L_{01} $ (in Proposition \ref{Construction of bitriple cover}). More precisely, we take the branch locus 
\begin{align*}
D_{01} \in \left| f_1\right|, D_{02} \in \left| f_1\right|, D_{22} \in \left| 3f_2\right|, D_{12} \in \left| 3f_3\right|
\end{align*}
\noindent
and $ D_{ij} = 0 $ for the other $ \left( i,j\right)  $ as Section \ref{The main construction 1 of surfaces with d = 3} except that now we require that the three divisors $ D_{01} $, $ D_{22} $, $ D_{12} $ meet transversally at a point $ P $. If we take a $ \mathbb{Z}_{3}^2 $-cover which is ramified on this branch locus, we obtain a surface with singularities. We resolve the singularities appearing in this case as follows: we blow up at the common point of $ D_{01} $, $ D_{22} $ and $ D_{12} $, the strict transforms of $ D_{\sigma} $ satisfy the conditions of Proposition \ref{Construction of bitriple cover} and Proposition \ref{Smoothness condition of bitriple cover}. Hence, the $ \mathbb{Z}_{3}^2 $-cover which is ramified on these strict transforms is smooth. The preimage of $ P $ is an elliptic curve with self-intersection $ -9 $.\\

Now we present precisely this variation where we obtain the surface described in the fifth row of Theorem \ref{THE MAIN THEOREM WITH D=3 PART 1}. Let $ Y_4 $ be a del Pezzo surface of degree $ 5 $. We consider the following smooth divisors:
\begin{align*}
D_{01}:= h_{14}, D_{02}:= f_{11}, D_{22}:= h_{24}+f_{21}+f_{22}, D_{12}:= h_{34}+f_{31}+f_{32}
\end{align*}
\noindent
and $ D_{ij} = 0 $ for the other $ \left( i,j\right)  $, where $ f_{11} \in \left| f_1\right| $, $ f_{21}, f_{22} \in \left| f_2\right| $, $ f_{31}, f_{32} \in \left| f_3\right| $ are distinct divisors of $ Y_4 $ such that no more than two of these divisors go through the same point and the following non-trivial divisors: 
$$	
\begin{tabular}{r r r r r}
$ L_{10} :=$ & $$&$f_2 $&$+2f_3  $&$ -e_4 $ \\
$ L_{01} :=$ & $f_1$&$+f_2 $&$+f_3 $&$ -e_4 $\\
$ L_{20} :=$ &  $$&$2f_2 $&$+f_3 $&$ -e_4 $\\
$ L_{02} :=$ & $f_1$&$+2f_2 $&$+2f_3 $&$ -2e_4 $\\		  
$ L_{11} :=$ &  $f_1$&$+2f_2 $&$ $&$ -e_4 $\\
$ L_{22} :=$ &  $f_1$&$+f_2 $&$ $&$ -e_4 $\\
$ L_{12} :=$ &  $f_1$&$ $&$+f_3 $&$ -e_4 $\\
$ L_{21} :=$ &  $f_1$&$ $&$+2f_3 $&$ -e_4 $.
\end{tabular} 
$$
\noindent
These divisors $ D_{\sigma} $, $ L_{\chi} $ satisfy the following relations:
$$
\begin{adjustbox}{max width=\textwidth}
\begin{tabular}{l l r l r l r l r r r r r r}
$ 3L_{10} $&$ \equiv $&&$   $& &$ $&$D_{22} $&$ +$&$2D_{12} $ &$ \equiv$ & $$&$3f_2 $&$+6f_3  $&$ -3e_4 $\\
$ 3L_{01} $&$ \equiv $ & $D_{01}$&$ +$ & $2D_{02}  $&$ +$&$D_{22} $&$ +$&$D_{12} $ & $ \equiv$ & $3f_1$&$+3f_2 $&$+3f_3 $&$ -3e_4 $\\
$ 3L_{20} $&$ \equiv $&&$ $&&$ $&$2D_{22} $&$ +$&$D_{12} $ & $ \equiv$ &  $$&$6f_2 $&$+3f_3 $&$ -3e_4 $\\
$ 3L_{02} $&$ \equiv$ & $2D_{01} $&$ +$ & $D_{02} $&$ +$&$2D_{22} $&$ +$&$2D_{12} $ & $ \equiv$ & $3f_1$&$+6f_2 $&$+6f_3 $&$ -6e_4 $\\
$ 3L_{11} $&$ \equiv $ & $D_{01}$&$ +$ & $2D_{02}  $  &$ +$&$2D_{22} $&& & $ \equiv$ &  $3f_1$&$+6f_2 $&$ $&$ -3e_4 $\\
$ 3L_{22} $&$ \equiv $ & $2D_{01}$&$ +$ & $D_{02}  $  &$ +$&$D_{22} $&& &$ \equiv$ &  $3f_1$&$+3f_2 $&$ $&$ -3e_4 $\\
$ 3L_{12} $&$ \equiv $ & $2D_{01}$&$ +$ & $D_{02}  $&$  $&&$ +$&$D_{12} $ &$ \equiv$ &  $3f_1$&$ $&$+3f_3 $&$ -3e_4 $\\
$ 3L_{21} $&$ \equiv $ & $D_{01}$&$ +$ & $2D_{02}  $&$  $&&$ +$&$2D_{12} $ &$ \equiv$ &  $3f_1$&$ $&$+6f_3 $&$ -3e_4 $.
\end{tabular}
\end{adjustbox}
$$
\noindent
Thus, by Proposition \ref{Construction of bitriple cover}, the divisors $ D_{\sigma} $, $ L_{\chi} $ define a $ \mathbb{Z}^2_3 $-cover $ \xymatrix{g: X \ar[r] & Y_4}  $. In addition, by Propositions \ref{Smoothness condition of bitriple cover} and \ref{invariants of bitriple cover}, the surface $ X $ is smooth and satisfies the following:
\begin{align*}
3K_X &\equiv g^*\left( h_{14} +h_{24}+h_{34}+ 2f_2 + 2f_3\right). 
\end{align*}		  
Since $ 3K_X $ is the pull-back of a nef and big divisor, the canonical divisor $ K_X $ is nef and big. Thus, the surface $ X $ is minimal and of general type. Furthermore, by Proposition \ref{invariants of bitriple cover}, the surface $ X $ possesses the following invariants:
\begin{align*}
K_X^2&= \left( h_{14} +h_{24}+h_{34}+ 2f_2 + 2f_3\right)^2 =21,\\
p_g\left( X\right) &=4, \chi\left( \mathcal{O}_X\right) =5, q\left( X\right)  = 0.   
\end{align*}

As in Section \ref{The main construction 1 of surfaces with d = 3}, because $ D_{01} + D_{22} + D_{12} $ has twelve double points and $ D_{02}\left( D_{01} + D_{22} + D_{12}\right) = 6 $, by Remark \ref{the type of singularities in the quotient surfaces} the quotient surface $ X_{1}:= X/\left\langle \left( 1,0\right) \right\rangle$ has twelve singular points of type $ A_2 $ and six singular points of type $ \frac{1}{3}\left( 1,1\right)  $. We consider the following maps
$$
\xymatrix{\overline{X}_1 \ar[r]^{\overline{g}_1} & \overline{Y}_4 \ar[r]^{\pi_1}& Y_4}
$$
where the map $ \pi_{1} $ is the blow-up at the intersection points between $ D_{02} $ and $ D_{22} +D_{12}$, the map $ \overline{g}_1 $ is the $ \mathbb{Z}_3 $-cover which is ramified on the strict transforms of $ D_{02}, D_{01}+D_{22}+D_{12} $. The surface $ \overline{X}_1 $ is a partial resolution of $ X_1 $ where the singularities of type $ \frac{1}{3}\left( 1,1\right)  $ are resolved. The surface $ \overline{X}_1 $ has twelve singular points of type $ A_2 $ and satisfies the following
\begin{align*}
3K_{\overline{X}_1} &\equiv \overline{g}_1^{*}\left( \pi_1^{*}\left( f_1 + 6l\right) -3e_2-3e_3-3e_4 -\sum_{i=1}^{6}{h_i}\right), 
\end{align*}
where $ h_1, h_2, \ldots, h_6 $ are the exceptional divisors of $ \pi_1 $. Because the surface $ \overline{X}_1 $ has only canonical singularities, the minimal resolution $ \tilde{X}_1 $ of $ X_1 $ has the following invariants:
\begin{align*}
K_{\tilde{X}_1}^2 &=K_{\overline{X}_1}^2 =\left( \pi_1^{*}\left( f_1 + 6l\right) -3e_2-3e_3-3e_4 -\sum_{i=1}^{6}{h_i}\right)^2= 5,\\
p_g\left( \tilde{X}_1\right) &=4, \chi\left( \mathcal{O}_{\tilde{X}_1}\right) =5, q\left( \tilde{X}_1\right)  = 0. 
\end{align*}		  
\noindent
Moreover, the surface $ \tilde{X}_1 $ satisfies the following
\begin{align*}
3K_{\tilde{X}_1} &\equiv \tilde{g}_1^{*}\left( \pi_1^{*}\left( f_1 + 6l\right) -3e_2-3e_3-3e_4 -\sum_{i=1}^{6}{h_i}\right), 
\end{align*}
\noindent
where $ \tilde{g}_1 $ is the composition of the map $ \overline{g}_1 $ and the resolution map which resolves the singular points of type $ A_2 $. Because the surface $ \tilde{X}_1 $ has no non-trivial $ 3 $-torsion, the canonical divisor of $ \tilde{X}_1 $ is
\begin{align*}
K_{\tilde{X}_1} &\equiv \theta +\tilde{g}_1^{*}\left( \pi_1^{*}\left( 2l\right) -e_2 - e_3-e_4\right),
\end{align*}
where $ \theta:=\tilde{g}_1^{*}\left( \pi_1^{*}\left( f_{11}\right)  -\sum\limits_{i=1}^{6}{h_i}\right)_{\text{red}} $. We remark that the divisor $ \pi_1^{*}\left( f_{11}\right)  -\sum\limits_{i=1}^{6}{h_i} $ is effective and is contained in the branch locus. So the divisor $ \theta $ is a $ \left( -2\right)  $-curve. We have that
\begin{align*}
h^{0}\left( \tilde{X}_1,\mathcal{O}_{\tilde{X}_1} \left( \tilde{g}_1^{*}\left( \pi_1^{*}\left( 2l\right) -e_2 - e_3-e_4\right) \right) \right) = 3.
\end{align*}

Because $ h^{0}\left( \tilde{X}_1,\mathcal{O}_{\tilde{X}_1} \left( \tilde{g}_1^{*}\left( \pi_1^{*}\left( 2l\right) -e_2 - e_3-e_4\right) \right) \right) \ne h^{0}\left( \tilde{X}_1,\mathcal{O}_{\tilde{X}_1} \left( K_{\tilde{X}_1} \right) \right) $, $ \theta K_{\tilde{X}_1} =0 $ and the linear system $ \left|  \tilde{g}_1^{*}\left( \pi_1^{*}\left( 2l\right) -e_2 - e_3-e_4\right) \right|  $ is base point free, the linear system $ \left|K_{\tilde{X}_1} \right|  $ is base point free. This implies that the canonical map $ \varphi_{\left| K_{\tilde{X}_1} \right| } $ is birational (see Section 1 \cite{MR1573789}). Thus the canonical map $ \varphi_{\left| K_{X} \right| } $ is a $ 3$-to-$1 $ map onto a quintic surface in $ \mathbb{P}^3 $ (see \cite[Theorem 4.6]{MR938031} and \cite[Theorem 1.5]{MR1076174}). 

Furthermore, because the $ \mathbb{Z}_3$-cover $ \xymatrix{X \ar[r]^{g_2} & X_1} $ is ramified on six singularities of type $ \frac{1}{3}\left( 1,1\right)  $ and twelve singularities of type $ A_2 $ and $ p_g\left( X\right)  = p_g\left( \tilde{X}_1\right) $,  by Proposition \ref{base-point-of-triple-cover}, the canonical system $ \left| K_X \right|  $ has six base points. So we obtain the surface described in the fifth row of Theorem \ref{THE MAIN THEOREM WITH D=3 PART 1}.

\begin{Remark}
	Similarly as above, if we require that the three branch components $ D_{01} $, $ D_{22} $, $ D_{12} $ meet transversally at one point and that the three branch components $ D_{02} $, $ D_{22} $, $ D_{12} $ meet transversally at one (two) point, we obtain the surface described in the sixth (the seventh) row of Theorem \ref{THE MAIN THEOREM WITH D=3 PART 1}. 
\end{Remark}

\section{Construction of the surface listed in Theorem \ref{THE MAIN THEOREM WITH D=3 PART 2}}		  
\subsection{Construction and computation of invariants}
Let $ l_1 $ be a line in $ \mathbb{P}^2 $, $ P_1 $ be a point which does not lie in the line $ l_1 $, and let $ P_2, P_3, P_4,P_5, P_6 $ be five distinct points in $ l_1 $. Let $ Y_6 $ be the blow-up of $ \mathbb{P}^2 $ at $ P_1, P_2, \ldots, P_6 $. We denote by 
\begin{center}
	\begin{tabular}{l l}
		$ l $: & the pullback of a general line in $ \mathbb{P}^2$,\\
		$ e_1, e_2, \ldots, e_6 $: & the exceptional divisors corresponding to $ P_1, P_2, \ldots, P_6 $, respectively,\\
		$ f_1, f_2, \ldots, f_6$: & the strict transforms of a general line through $ P_1, P_2, \ldots, P_6 $, respectively,\\
		$ h_{23456} $: & the strict transform of the line $ l_1 $,\\
		$ h_{12}, h_{13}, \ldots, h_{16} $: & the strict transforms of the lines $ P_1P_2, P_1P_3, \ldots, P_1P_6 $, respectively.
	\end{tabular}
\end{center}

We consider the following smooth divisors:
\begin{align*}
D_{01}:= h_{14}+h_{15}+h_{16}, D_{02}:= h_{23456}, D_{22}:= f_{21}+f_{22}+e_{3}, D_{12}:= f_{31}+f_{32}+e_{2}
\end{align*}
\noindent
and $ D_{ij} = 0 $ for the other $ \left( i,j\right)  $, where $ f_{21}, f_{22}\in \left| f_2\right| $, $ f_{31}, f_{32}\in \left| f_3\right| $ are distinct divisors in $ Y_6 $ such that no more than two of these divisors go through the same point and the following non-trivial divisors: 
$$	
\begin{tabular}{r r r r l}
$ L_{10} :=$ & $$&$$&$f_3  $&$ +l $ \\
$ L_{01} :=$ & $f_1$&$+f_2 $&$+f_3 $&$ -e_4 -e_5 -e_6$\\
$ L_{20} :=$ &  $$&$f_2 $&$ $&$ +l $\\
$ L_{02} :=$ & $2f_1$&$+f_2 $&$+f_3 $&$+f_4-e_5 -e_6 $\\		  
$ L_{11} :=$ &  $f_1$&$+2f_2 $&$ $&$ -e_4 -e_5 -e_6$\\
$ L_{22} :=$ &  $2f_1$&$+f_2 $&$ $&$ -e_4 -e_5 -e_6$\\
$ L_{12} :=$ &  $2f_1$&$ $&$+f_3 $&$ -e_4 -e_5 -e_6$\\
$ L_{21} :=$ &  $f_1$&$ $&$+2f_3 $&$ -e_4 -e_5 -e_6$.
\end{tabular} 
$$
\noindent
These divisors $ D_{\sigma} $, $ L_{\chi} $ satisfy the following relations:
$$
\begin{adjustbox}{max width=\textwidth}
\begin{tabular}{l l r l r l r l r r r r r l}
	$ 3L_{10} $&$ \equiv $&&$   $& &$ $&$D_{22} $&$ +$&$2D_{12} $ &$ \equiv$ & $$&$$&$3f_3  $&$ +3l $\\
	$ 3L_{01} $&$ \equiv $ & $D_{01}$&$ +$ & $2D_{02}  $&$ +$&$D_{22} $&$ +$&$D_{12} $ &$ \equiv$ & $3f_1$&$+3f_2 $&$+3f_3 $&$ -3e_4 -3e_5 -3e_6$\\
	$ 3L_{20} $&$ \equiv $&&$ $&&$ $&$2D_{22} $&$ +$&$D_{12} $ &$ \equiv$ &  $$&$3f_2 $&$ $&$ +3l $\\
	$ 3L_{02} $&$ \equiv$ & $2D_{01} $&$ +$ & $D_{02} $&$ +$&$2D_{22} $&$ +$&$2D_{12} $ &$ \equiv$ & $6f_1$&$+3f_2 $&$+3f_3 $&$+3f_4-3e_5 -3e_6 $\\
	$ 3L_{11} $&$ \equiv $ & $D_{01}$&$ +$ & $2D_{02}  $  &$ +$&$2D_{22} $&& &$ \equiv$ &  $3f_1$&$+6f_2 $&$ $&$ -3e_4 -3e_5 -3e_6$\\
	$ 3L_{22} $&$ \equiv $ & $2D_{01}$&$ +$ & $D_{02}  $  &$ +$&$D_{22} $&& &$ \equiv$ &  $6f_1$&$+3f_2 $&$ $&$ -3e_4 -3e_5 -3e_6$\\
	$ 3L_{12} $&$ \equiv $ & $2D_{01}$&$ +$ & $D_{02}  $&$  $&&$ +$&$D_{12} $ &$ \equiv$ &  $6f_1$&$ $&$+3f_3 $&$ -3e_4 -3e_5 -3e_6$\\
	$ 3L_{21} $&$ \equiv $ & $D_{01}$&$ +$ & $2D_{02}  $&$  $&&$ +$&$2D_{12} $ &$ \equiv$ &  $3f_1$&$ $&$+6f_3 $&$ -3e_4 -3e_5 -3e_6$.
	\end{tabular}
	\end{adjustbox}
	$$
	\noindent
	Thus, by Proposition \ref{Construction of bitriple cover}, the divisors $ D_{\sigma} $, $ L_{\chi} $ define a $ \mathbb{Z}^2_3 $-cover $ \xymatrix{g: X \ar[r] & Y_6}  $. Moreover, by Propositions \ref{Smoothness condition of bitriple cover} and \ref{invariants of bitriple cover}, the surface $ X $ is smooth and satisfies the following:
	\begin{align*}
	3K_X &\equiv g^*\left( 3l+3f_1 + h_{23456}\right). 
	\end{align*}		  
	Since $ 3K_X $ is the pull-back of a nef and big divisor, the canonical divisor $ K_X $ is nef and big. Thus, the surface $ X $ is minimal and of general type. Furthermore, by Proposition \ref{invariants of bitriple cover}, the surface $ X $ possesses the following invariants:
	\begin{align*}
	K_X^2&=\left( 3l+3f_1 + h_{23456}\right)^2= 35,\\
	p_g\left( X\right) &=6, \chi\left( \mathcal{O}_X\right) =7, q\left( X\right)  = 0.   
	\end{align*}
	
	\subsection{The canonical map and the canonical image}
	Let $ X $ be as above and let $ \Sigma $ denote the image of the canonical map $ \varphi_{\left| K_X \right| } $. We show that the degree $d$ of the canonical  map is $ 3 $ and that the canonical image $ \Sigma $ is a surface of general type. We consider the cyclic subgroup $ \Gamma:= \left\langle \left( 1,0\right) \right\rangle  $ of $ \mathbb{Z}_3^2 $. Because $ h^0\left( L_{ \chi} + K_{Y_6} \right)  = 0 $ for all $ \chi \notin \Gamma^\perp $, by Proposition \ref{Factorization of the canonical map through a quotient surface}, the canonical map $ \varphi_{\left| K_X \right| } $ is essentially the composition of the quotient map $g_{2}: \xymatrix{X \ar[r]& X_{1}:= X/\left\langle \left( 1,0\right) \right\rangle } $ with the canonical map $\varphi_{\left| K_{\tilde{X}_{1}} \right| }$ of the minimal resolution $ \tilde{X}_{1} $ of $ X_1 $.
	
	The surface $ X_1 $ is the $ \mathbb{Z}_3 $-cover $ \xymatrix{g_1 : X_1 \ar[r]& Y_6} $ with the following building data:
	\begin{align*}
	D_{1}&:= D_{01} + D_{22} +D_{12}, &D_{2}&: = D_{02}\\		  
	L_{1}&: = L_{01}, &L_{2}&: = L_{02}.		  
	\end{align*}
	Because the divisor $ D_1 $ has twenty double points and $ D_1D_2 = 2 $, by Remark \ref{the type of singularities in the quotient surfaces} the quotient surface $ X_1 $ has twenty singular points of type $ A_2 $ and two singular points of type $ \frac{1}{3}\left( 1,1\right)  $. We consider the following maps
	$$
	\xymatrix{\overline{X}_1 \ar[r]^{\overline{g}_1} & \overline{Y}_6 \ar[r]^{\pi_1}& Y_6}
	$$
	where the map $ \pi_{1} $ is the blow-up at the intersection points between $ D_{02} $ and $ D_{22} +D_{12}$, the map $ \overline{g}_1 $ is the $ \mathbb{Z}_3 $-cover which is ramified on the strict transforms of $ D_{02}, D_{01}+D_{22}+D_{12} $. The surface $ \overline{X}_1 $ is a partial resolution of $ X_1 $ where the $ \frac{1}{3}\left( 1,1\right)  $ singularities are resolved. The surface $ \overline{X}_1 $ has twenty singular points of type $ A_2 $ and satisfies the following
	\begin{align*}
	3K_{\overline{X}_1} &\equiv \overline{g}_1^{*}\left( \pi_1^{*}\left( 3l+3f_1 + h_{23456}\right)  -h_1-h_2\right), 
	\end{align*}
	where $ h_1, h_2$ are the exceptional divisors of $ \pi_1 $. Because the surface $ \overline{X}_1 $ has only canonical singularities, the minimal resolution $ \tilde{X}_1 $ of $ X_1 $ has the following invariants:
	\begin{align*}
	K_{\tilde{X}_1}^2 &=K_{\overline{X}_1}^2=\left( \pi_1^{*}\left( 3l+3f_1 + h_{23456}\right)  -h_1-h_2\right)^2= 11,\\
	p_g\left( \tilde{X}_1\right) &=6, \chi\left( \mathcal{O}_{\tilde{X}_1}\right) =7, q\left( \tilde{X}_1\right)  = 0. 
	\end{align*}		  
	\noindent
	Moreover, the surface $ \tilde{X}_1 $ satisfies the following
	\begin{align*}
	3K_{\tilde{X}_1} &\equiv \tilde{g}_1^{*}\left( \pi_1^{*}\left( 3l+3f_1 + h_{23456}\right)  -h_1-h_2\right), 
	\end{align*}
	\noindent
	where $ \tilde{g}_1 $ is the composition of the map $ \overline{g}_1 $ and the resolution map which resolves the singular points of type $ A_2 $.
	
	Analogously to Section \ref{The main construction 1 of surfaces with d = 3}, the surface $ \tilde{X}_1 $ has no non-trivial $ 3$-torsion. So we get
	\begin{align*}
	K_{\tilde{X}_1} &\equiv \theta +\tilde{g}_1^{*}\left( \pi_1^{*}\left( l\right) \right) + \tilde{g}_1^{*}\left( \pi_1^{*}\left( f_1\right) \right),
	\end{align*}
	where $ \theta:=\tilde{g}_1^{*}\left( \pi_1^{*}\left( h_{23456}\right)  -h_1-h_2\right)_{\text{red}} $. We remark that the divisor $ \pi_1^{*}\left( h_{23456}\right)  -h_1-h_2 $ is effective and is contained in the branch locus. So the divisor $ \theta $ is a $ \left( -2\right)  $-curve. We have that
	\begin{align*}
	h^{0}\left( \tilde{X}_1,\mathcal{O}_{\tilde{X}_1} \left( \tilde{g}_1^{*}\left( \pi_1^{*}\left( l\right) \right) + \tilde{g}_1^{*}\left( \pi_1^{*}\left( f_1\right) \right) \right) \right) = 5.
	\end{align*}

	Because $ h^{0}\left( \tilde{X}_1,\mathcal{O}_{\tilde{X}_1} \left( \tilde{g}_1^{*}\left( \pi_1^{*}\left( l\right) \right) + \tilde{g}_1^{*}\left( \pi_1^{*}\left( f_1\right) \right) \right) \right) \ne h^{0}\left( \tilde{X}_1,\mathcal{O}_{\tilde{X}_1} \left( K_{\tilde{X}_1} \right) \right) $, $ \theta K_{\tilde{X}_1} =0 $ and the linear system $ \left|  \tilde{g}_1^{*}\left( \pi_1^{*}\left( l\right) \right) + \tilde{g}_1^{*}\left( \pi_1^{*}\left( f_1\right) \right) \right|  $ is base point free, the linear system $ \left|K_{\tilde{X}_1} \right|  $ is base point free. Because the linear system $ \left|K_{\tilde{X}_1} \right|  $ is base point free, we get the following equality 
	\begin{align*}
	\degree\left( \varphi_{\left| K_{\tilde{X}_1} \right|}\right) \degree\left( \Sigma \right) = K_{\tilde{X}_1}^2 =11.
	\end{align*}
	\noindent
	This implies that the map $ \varphi_{\left| K_{\tilde{X}_1} \right| } $ is birational because $ K_{\tilde{X}_1}^2 $ is a prime number. In fact, since the image $ \Sigma $ is contained in $ \mathbb{P}^5 $, $ \degree\left( \Sigma \right) \ge 4  $. So $ \degree\left( \varphi_{\left| K_{\tilde{X}_1} \right|}\right) = 1 $. In addition, because the linear system $ \left|K_{\tilde{X}_1} \right|  $ is base point free, the canonical map $ \varphi_{\left| K_{\tilde{X}_1} \right| } $ is a morphism.

	Finally, because the $ \mathbb{Z}_3$-cover $ \xymatrix{X \ar[r]^{g_2} & X_1} $ is ramified on the $ \mathbb{Z}_3 $-set of two singularities of type $ \frac{1}{3}\left( 1,1\right)  $ and twenty singularities of type $ A_2 $, and $ p_g\left( X\right)  = p_g\left( \tilde{X}_1\right) $,  by Proposition \ref{base-point-of-triple-cover}, the canonical system $ \left| K_X \right|  $ has two base points. So we obtain the surface described in Theorem \ref{THE MAIN THEOREM WITH D=3 PART 2}.
	
	\begin{Remark}
		In all the constructions above, the base points come from the singularities of type $ \frac{1}{3}\left( 1,1\right)  $ of the quotient surface $ X_1 $ (see Proposition \ref{base-point-of-triple-cover}).
	\end{Remark}

	

\Addresses

\Addressesb

\Addressesc

\end{document}